# The Trade-off Strategy between Railroads and Customers: A Case Study for Low-frequency Entire Trains

Boliang Lin*, Zexi Zhang

*School of Traffic and Transportation, Beijing Jiaotong University, Beijing 100044, China*

**Abstract:** Some large freight railroads ship a number of shipments over the rail network annually. To reduce unnecessary reclassifications of shipments on their routes, each railroad is willing to operate the entire train for an individual shipment. In other words, the motivation for providing the entire train service lies in a simple realization that "door to door" transportation (directly from origin to destination) can reduce operating costs by decreasing classification. However, this mode will increase inventory costs for customers when commodities are transported by low-frequency entire train services. Thus, this study proposes the trade-off strategy to keep a balance between saving operating costs of railroads and increasing inventory costs of customers. We analyze the revenue and losses after a shipment shifting from the transfer transportation which contains a series of train services to the direct transportation by entire train service.

## 1. Introduction

The entire train service is an important mode to transport bulk commodities directly from origin to destination without any intermediate handling at classification yard, and usually contains a single bulk commodity (i.e., coal, iron ore and grain) and a shipment for a single customer (Ravindra K. Ahuja 2014). This could bring benefits to railroads by reducing reclassification costs. Some large railroads usually provide entire train services to achieve "door to door" transportation for some shipments with large volumes of freight traffic. However, those remaining shipments will be shipped by transfer transportation and have to be reclassified at marshalling yards on their routes. In this case, operating costs will usually increase. In fact, Keaton (1991) has found that railroads and customers have accepted the inefficiencies of classification yards as the price of shipment consolidation.

In general, the entire train will be provided when size of the shipment is larger than one train per day. In some cases, it is also beneficial for shipment of small volume to provide the entire train below one train per day, and that can be called low-frequency entire train. The low-frequency entire train differs from the entire train in that this transportation mode is provided with a longer service time interval. Compared with common transfer transportation for small volumes of freight traffic, this shipping strategy can reduce operating costs by decreasing classification.

Based on saving classification costs, railroads could adopt a flexible pricing policy to entice more customers to choose the entire train. The operating revenue of a railroad must be sufficient to cover the railroad's operating costs and provide an acceptable return on infrastructure investment. Therefore, this paper innovatively designs a pricing discount strategy for the low-frequency entire train. Namely, railroads provide train services at a discount of rail charges to customers achieving win-win situation between railroads and customers.

* Corresponding author at: School of Traffic and Transportation, Beijing Jiaotong University, Beijing 100044, China. Tel./fax: +86 10 51687149. E-mail address: bllin@bjtu.edu.cn

It is usually very complicated to give a suitable discount that satisfies profitable for railroads and eventually attract more customers to choose entire train services. Although operating the entire train service is beneficial to railroads, this will significantly increase inventory costs for customers. That is, when the demand of customers per day is much fewer than the size of a train, the excess freight need to be stored. Therefore, the discount of rail charges should be adopted to compensate inventory costs of customers.

The objective of this paper is to aid railroads in determining the discount pricing considering interests with customers. The contributions of this study are summarized as follows. Firstly, we propose the low-frequency entire train service adopting a discount pricing strategy. Secondly, this paper discusses four trade-off strategies for keeping a balance between the saving operating costs of railroads and the increasing inventory costs of customers.

The rest of this paper is as follows: Section 2 presents a brief literature survey on discount pricing strategies, the shipping strategies and theory of inventory. In Section 3, the paper describes theories of the shipping strategies, rail charges and inventory. In Section 4, several typical trade-off strategies are discussed for balancing the increasing inventory costs of the customers and saving reclassification costs by operating the low-frequency entire train services. Finally, Section 5 concludes and outlines the future research prospect.

## 2. Literature review

Several researches related to this paper can be divided into three different sub-sections: discount pricing strategies, the shipping strategies and theory of inventory.

### 2.1 Discount pricing strategies

The literature which closely relates to discount pricing strategies should be mentioned. Discount pricing strategies have been widely adopted by suppliers to attract more customers to order in large batches. Among them, there are lots of research on quantity discount scheme. Crowther (1964) first proposed a quantity discount scheme simultaneously increasing a supplier's profit and reducing buyer's costs. Monahan (1984) was one of the earliest to adjust pricing schedule to encourage major customers to increase their order size, and optimize a factor to maximize supplier's net profit. Most subsequent research (Goyal, 1987. Banerjee,1986) also generalized Monahan's model by analyzing supplier's carrying costs. Gradually, the discount pricing strategies also began to shift from solely considering vendor's pricing schedule to comprehensive pricing schedule incorporating vendor's inventory carrying costs. Shinn et al. (1996) also considered a freight cost which is a function of the lot-size, and focused on the problem of simultaneously determining the retailer's lot size and optimal price that contains a fixed set-up cost and a freight cost offered quantity discount.

Pricing discount strategies are used not only by vendors but also by common carriers to attract more customers. Tsao and Sheen (2012) discussed the possibility of realizing considerable savings by considering a joint multi-item replenishment policy and weight freight cost discounts from both individual and channel perspectives. Li et al. (2012) established mixed integer programming formulations with quantity discounts for the airfreight forwarder's shipment planning problem. Nguyen et al. (2014) proposed a purchase incentive to increase buyers' order sizes and using multi-stop routes to promote truck utilization. Different from the widely studied in trucking and airline shipments, Qiu et al. (2019) studied a rail transportation pricing problem with all-units quantity discounts, and formulated a Stackelberg game model concerning one dry port and multiple shippers.

Compared with the quantity discount strategy, this research focuses on the discount pricing strategies for the choice of modes in rail transportation. The American railroad has been deregulated since the Staggers Rail Act of 1980. And a great deal of subsequent research has focused on the change in rail rates after the reform and deregulation. Burton (1993) and Wilson (1994) summarized that railroad deregulations did promote the falling rail rates differing from general feeling that this change would lead to higher railroad rates. And they also found that shippers would benefit from decreasing rail freight rates and changing the characteristics of their shipments. Faced with the decline in railroad rates, Dennis (2001) analyzed the importance of each of these factors. Ellig (2002) focused on unreasonable rates for captive shippers and finally proposed the method of cost reduction. Based on railroads deregulation, more literature develops new rail rates pricing strategies. Sparger and Prater (2012) mentioned that unit trains and shuttle trains adopt new rail rate indices using grain as an example, and the innovation indices include a higher level of details in contrast to other annualized indices. Ndembe (2015) discussed railroad pricing behavior and found that intermodal competition and shuttle trains can reduce rate over time. Similar to Sparger and Prater, Hyland et al. (2016) established conceptual and mathematical models about grain transportation supply considering travel time, cost and capacity, and demonstrated benefits of shuttle train comparing the conventional service.

2.2 The shipping strategies

As Sparger and Prater (2012) mentioned previously, railroads began adopting pricing incentives for larger car shipments after the 1980s. To further reduce transportation cost and increase efficiency of commodities transportation, Frittelli (2005), Sparger and Prater (2012) and Ravindra K. Ahuja (2014) et.al defined the transportation modes of bulk commodities (including the term unit train, the term shuttle train and the term entire train). Generally speaking, unit trains refer to a block of 50-54 railcar shipments which obtain a lower per-railcar pricing while shuttle trains refer to an entire train of 90-120 railcar shipments which consisting of the locomotives and crew ( Frittelli 2005, Sparger and Prater 2012). Similar to the shuttle train, the entire trains are provided for a single shipment containing a single commodity from an individual origin and destination rather than multiple origins and destinations (Ravindra K. Ahuja 2014).

As the shipment sizes are not sufficient for filling entire trains, they will be shipped by transfer transportation and have to be intermediate handled at classification yards on their routes. A large volume of literature which is closely related to these shipping strategies should be noted as well. The freight transfer transportation problem in railway is also known as the railroad blocking problem. The literature on railroad blocking problem is abundant and rich. Railroad blocking is critical for railroads that minimize total transportation and handling cost (Ravindra K. Ahuja et al. 2007). To prevent shipments from being excessively reclassified at every classification yard they pass through, railroads group a set of railcars with the same origin and destination at a classification yard to create origin-destination pairs, known as blocks. However, the origin and the destination of each individual shipment may differ from that of blocks. Therefore, these shipments which are different from the destination of the block will be reclassified. Ravindra K. Ahuja et al. (2007)

developed large-scale neighborhood search to solve real-life railroad blocking problems including a variety of practical and business constraints. Lin et al. (2012) proposed a model for the freight train connection services problem of a large-scale China rail system. Finally, they applied the simulated annealing algorithm to solve the problem and obtained decent result. Zhu et al. (2014) addressed the scheduled service network design problem for freight rail transportation based on a cyclic three-layer time-space network. Chen et al. (2018) established an exact linear binary model to minimize the total sum of train accumulation cost and car classification cost. And they proposed a tree-based decomposition algorithm, which can derive high-quality solutions. In practice, classification involves this process of sorting railcars and combining shipments. For the commodities by railway, the shipment may pass through many classification yards to be sorted and grouped. Boysen et al. (2012) concluded the previous literature about operational processes at classification yards during the past 40 years, which will help us to understand the classification process better.

2.3 Inventory management

The analysis of inventory costs for customers is still essential for our research. The additional increasing inventory costs for customers decide whether they choose the entire train service. After all, the core is to keep a balance between saving the operating costs of railroads and the increasing inventory costs of customers.

The existing literature on inventory theory is extensive and focuses particularly on the inventory model. Baumol et al. (1970) first systematically studied the integration of transportation and inventory costs and presented the inventory-theoretic models. A relatively large volume of literature used their inventory-theoretic approach as a basis for further development. Carried out by Swenseth and Godfrey (2002) formulated a realistic transportation costs function and proved that the straightforward freight rate function could be incorporated into inventory replenishment decisions without compromising the accuracy of the decision. Silver et al. (2008) had an overview of inventory management and listed contributions of fundamental and applied theory within the dimensions. Baller et al. (2019) incorporated approximated transportation costs by using classical schemes into a Dynamic-Demand Joint Replenishment Problem (DJRP) and assumed that the supplier paid a fixed fee for replenishing a customer. They analyzed the data from test instances and concluded that using a Dynamic-Demand Joint Replenishment Problem cost structure can lead to greater cost savings. Other studies paid attention to adding transportation costs into inventory costs. In this way, the research to analyze a systematic optimization of transportation and the related enterprise inventory can be seen in Hill and Omar (2006) and Glock (2012).

3. Problem description

3.1 The entire train and the transfer transportation

There are two different shipping strategies for delivering commodities from their origins to their respective destinations. Strategy 1, the railroad will operate an entire train service (provided between loading and unloading stations, see green line of Strategy 1 in Fig.1) without reclassification at any classification yard they pass through.

Strategy 2, transfer transportation strategies for the remaining commodities are as follows: They will be firstly sent to the closest classification yard $Y_1$ by a local train and be grouped together with others to form a block according to their destinations. A block is associated with an origin-destination pair that will be assigned to appropriate train services. There are also many different combinations for shipping these remaining commodities in transfer transportation strategies. For example: Strategy 2A, they may be delivered to the classification yard $Y_{k+1}$ by a series of direct train services (dispatched between a pair of yards bypassing several classification yards, see orange lines of Strategy 2A in Fig.1); or else, they may be sent by link train services (provided between the adjacent classification yards, see purple lines of Strategy 2B in Fig.1) and direct train service 3 until reaching classification yard $Y_{k+1}$. And then, the shipments will be delivered to the unloading station by a local train. In general, local train services are also provided between the adjacent classification yards in addition to link train services. In practical, if some bulk commodities satisfy the sufficient condition (i.e., more than 500 thousand tons per year), the entire train service will be provided for them without optimization. The remaining commodities will usually be delivered by Strategy 2A and 2B.

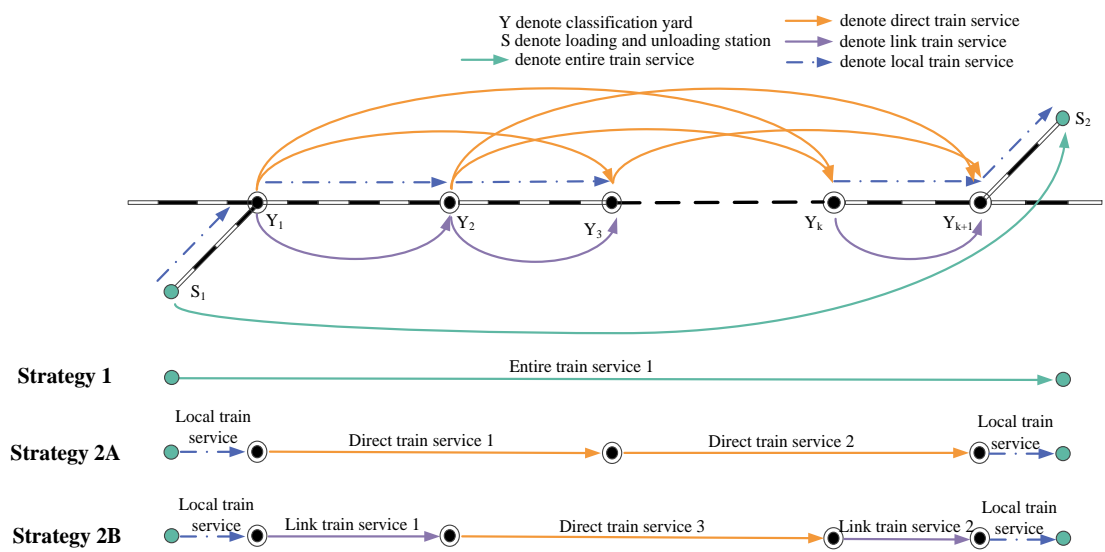

Figure 1. The delivering process for one shipment

Among these basic strategies above, the entire train service can obviously reduce the operating costs and save intermediate delay time by decreasing reclassification. Thus, railroads are willing to operate the entire train for each individual shipment. Clearly, it is beneficial for railroads for large volumes of freight traffic to provide the entire train. Sometimes, when daily demand of customers may be too small, these commodities are still transported by the entire train service. Although railroads decrease reclassification cost, inventory costs of customers usually increase compared with common transfer transportation. Generally speaking, entire trains are usually provided more than one train per day. Therefore, we propose the concept of the low-frequency entire train service. That is, the entire train with a longer service time interval (usually 2 days or even longer) not only allows railroads reduce reclassification cost but also meets daily demand of customers with a relatively low average inventory.

To have a better intuitive understanding of difference in the entire train and transfer transportation, we further analyze the component of railway transportation cost. The total cost of railway transportation can be divided into four parts: the loading cost, the car-miles cost, the intermediate handling cost (or the reclassification cost) and the unloading cost. Among them, the car-miles cost denotes the cost of a shipment sent from its origin to its destination on its route (including car hire, locomotive and crew cost). And the intermediate handling cost denotes the cost of the reclassification of a shipment at a classification yard. Comparing the entire train service and the transfer transportation, the car-miles cost is mostly equal, but the intermediate handling cost will be saved when the shipment is delivered by the entire train service. Loading and unloading costs increasing make up a miniscule portion comparing with saving reclassification costs. Hence, it is not influential for railroads to operate the entire train.

In these classification yards, the incoming trains, which may contain a great deal of individual shipments, are sorted and grouped. That is, inbound trains are disassembled and outbound trains are generated. Therefore, these classification processes delay the movement of shipments and are rather time and cost consuming. In this case, shipments delivered by the entire train services could bring benefits to railroads by reducing reclassification. As Ravindra K. Ahuja mentioned, the reclassification costs are different at different classification yards. In general, larger classification yards have lower reclassification costs and smaller classification yards have higher reclassification costs. For a typical U.S. railroad, the average reclassification cost is $40 per car, and a shipment with an average distance of 500 miles from its origin to its destination may be reclassified 2.5 times. Assuming that there is a shipment per day containing about 15 railcars, we could save about $547.5 thousand in reclassification cost for a shipment per year. These potential savings encourage railroads to provide the entire train for shipment of small volume. By operating the entire train without reclassification, the railroad can not only reduce the transportation costs but also save the shipping time. In order to synthetically consider the time and costs saved in the process of operating the entire train, we will transform from the saved time to the reduced costs by means of a conversion rate between time and cost.

Based on saving transportation cost and time, railroads could offer pricing incentive to entice customers to choose the entire train. In general, the operating revenue of the railroad must be sufficient to cover the railroad's operating costs, as well as provide an acceptable return on infrastructure investment. Thus, it is possible to provide the low-frequency entire train service adopting a discount pricing strategy.

## 3.2 Inventory costs and the low-frequency entire train

As mentioned above, from the perspective of transportation costs of railroads, not only shipments of large volume of freight traffic benefit from operating the entire train, shipments of small volume also benefit from reducing reclassification cost. However, operating the entire train service will significantly increase inventory costs of customers comparing with the transfer transportation. Typically, railroads provide transfer transportation services to customers with small daily demand and do not cause short-term surges for inventory of customers. In some cases, customers with small daily demand are also willing to choose the entire train below one train per day when inventory costs are compensated. In this section, we will analyze the change of inventory costs for customers due to the shift of train service mode.

As shown in Figure 2A, the entire train is provided one train per day, and the low-frequency entire train is provided one train per three days. Compared with the entire train, the low-frequency entire train could bring acceptable inventory cost increasing for customers. Let $Q^{\text{Train}}$ be the rated load of an entire train. The difference between the low-frequency entire train and the entire train is that the low-frequency entire train has a longer service time interval. Thus, the rated load of the entire train and the low-frequency entire train is equal.

Then $Q^{\text{Train}}$ is defined.

$$Q^{\text{Train}} = Q^{\text{Car}} m_{ij} \tag{1}$$

Where $m_{ij}$ is the size of the train dispatched from $i$ to $j$. And $Q^{\text{Car}}$ is the rated load of a unit railcar.

In practice, train services are evenly provided to customers with a fixed frequency (see Figure 2B). Let $Q^{\text{The non-direct train}}$ be the volume of freight traffic transported by the non-entire train every time, and $Q^{\text{Daily demand}}$ denotes the volume of daily demand of customers per day (Generally speaking, $Q^{\text{The non-direct train}} \geq Q^{\text{Daily demand}}$). As the Figure 2B shown, the volume of freight traffic transported by the non-entire train meets the volume of demand of customers during the contract's period. Obviously, inventory costs reach minimal in an extreme situation when commodities are delivered on daily demand of customers (when $Q^{\text{The non-direct train}} = Q^{\text{Daily demand}}$). Therefore, we assume that commodities are delivered by transfer transportation at minimal inventory costs. And $Q^{\text{Daily demand}}$ can be stated as follows.

$$Q^{\text{Daily demand}} = Q^{\text{Car}} N_{ij} \tag{2}$$

Where $N_{ij}$ denotes the number of daily railcars which origin at yard $i$ and are destined to yard $j$.

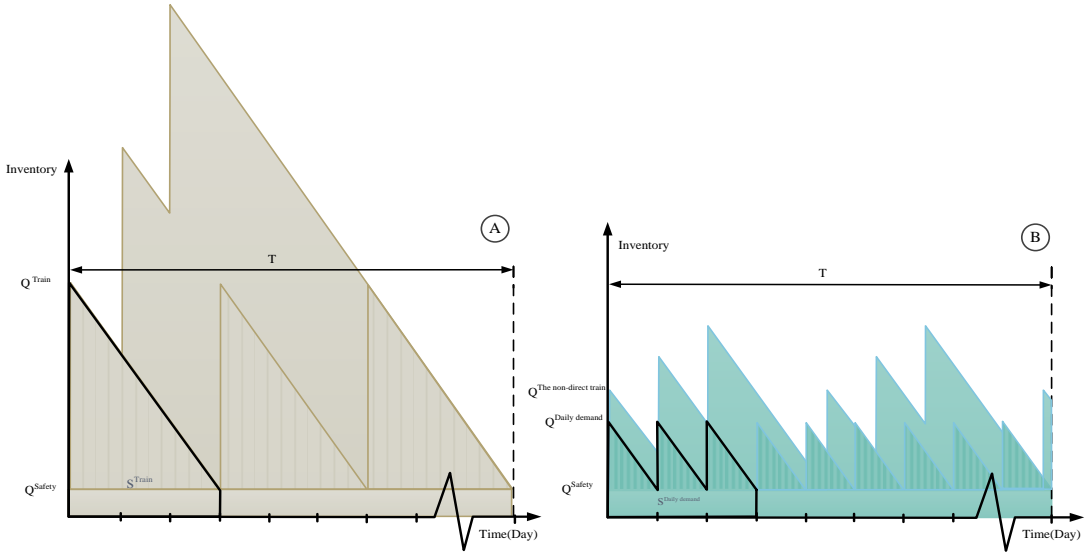

Figure 2. (A) The stock transported by the low-frequency entire train.

(B) The stock transported by the non-entire train.

Assuming that the loading and unloading time at the loading and unloading station is very short, that is, the consumption curve at the unloading and loading is linearly continuous (see Figure 2). The horizontal coordinate indicates the time, and the vertical coordinate represents the inventory. The total stock is the integral of the inventory curve. Because of the specificity of the inventory curve, the integral of the inventory curve can be approximated by the area of the figure (e.g., $S^{\text{Daily demand}}$ and $S^{\text{Train}}$ in Fig.2). The area of the figure denotes the stock during the same contract's period.

Let $Q^{\text{Safety}}$ represent the safety stock, and we assume that the length of a period is equal to three days. Then, $S^{Train}$ is the stock transported by the low-frequency entire train. This can be defined by,

$$S^{Train} = \frac{3}{2}Q^{Train} + 3Q^{Safety}$$
$$= \frac{3}{2}Q^{Car}m_{ij} + 3Q^{Safety} \quad (3)$$

Similarly, the $S^{Daily\ demand}$ is the stock transported by the non-entire train. This can be stated as,

$$S^{Daily\ demand} = \frac{3}{2}Q^{Daily\ demand} + 3Q^{Safety}$$
$$= \frac{3}{2}Q^{Car}N_{ij} + 3Q^{Safety} \quad (4)$$

Therefore, $\Delta S$ is the extra increasing stock per day by customers choosing the low-frequency entire trains, and this can be stated as,

$$\Delta S = \frac{1}{3}\left(S^{Train} - S^{Daily\ demand}\right) = \frac{1}{2}Q^{Car}(m_{ij} - N_{ij}) \quad (5)$$

This section analysis that when the customers choose the low-frequency entire train services, customers will have to pay more extra inventory costs. On the one hand, assuming that inventory costs of customers are minimal until changing the train service mode, this will ensure that railroads can satisfy more customers by adopting pricing incentives. On the other hand, we propose that the entire train with a longer service time interval (the low-frequency entire train) could bring acceptable inventory cost increasing for customers.

## 4. Methodology

### 4.1 Parameters

This section quantitatively analyzes whether it is profitable to operate low-frequency entire trains from the perspective of the railroads and the customers. We assume that daily demand and the volume of freight traffic annually within a given contract period is constant regardless of whether the low-frequency entire train is operated.

Railroads may save the reclassification time and the intermediate handling costs by decreasing reclassification. However, operating the entire train service will significantly increase inventory costs of customers. Thus, we propose the trade-off strategy in this section to calculate a suitable discount that satisfies profitable for railroads and attract more customers to choose low-frequency entire train services.

The relevant notations are described as follows.

| Sets | |
|---|---|
| $V$ | The set of classification yards in a rail network |
| Parameters | |
| $\Delta H$ | The difference between the cost-saving and discounts when the railroad providing a low-frequency entire train with pricing incentive, in yuan. |
| $\Delta E$ | The difference of cost and time whether the railroad provides a low- |
| | frequency entire train for a customer, in yuan. |
| $\beta$ | The discount parameter when the customers choose the low-frequency entire train services. |
| $P$ | The rail charges, in yuan. |
| $\Delta E^{\text{loading}}$ | The difference in loading cost at the loading station whether the railroad provides a low-frequency entire train for a customer, in yuan. |
| $\Delta E^{\text{unloading}}$ | The difference in unloading cost at the unloading station whether the railroad provides a low-frequency entire train for a customer, in yuan. |
| $\Delta E^{\text{reclassification}}$ | The difference in synthetically considering reclassification cost and time at the classification yards which they pass through whether the railroad provides a low-frequency entire train for a customer, in yuan. |
| $\Delta E^{\text{car-miles}}$ | The difference in car-miles cost of transporting a shipment from its origin to its destination whether the railroad provides a low-frequency entire train for a customer, in yuan. |
| $\Delta G^{\text{loading}}$ | The difference in loading time at the loading station whether the railroad provides a low-frequency entire train for a customer, in car-hour. |
| $\gamma$ | Conversion parameter between time and cost, in yuan/car-hour. |
| $C^{\text{loading}}$ | The extra cost of loading for operating a low-frequency entire train, in yuan/car. |
| $\Delta G^{\text{unloading}}$ | The difference in unloading time at the unloading station whether the railroad provides a low-frequency entire train for a customer, in car-hour. |
| $C^{\text{unloading}}$ | The extra cost of unloading for operating a low-frequency entire train, in yuan/car. |
| $\Delta G^{\text{reclassification}}$ | The difference in reclassification time at the classification yards which they pass through whether the railroad provides a low-frequency entire train for a customer, in car-hour. |
| $\Delta C^{\text{reclassification}}$ | The difference in the single cost at the classification yards which they pass whether the railroad provides a low-frequency entire train for a customer, in yuan. |
| $q$ | The number of classification yards that should have been reclassified according to the train formation plan from the classification yard, one. |
| $t_k^{\text{broken up}}$ | The parameter of the train broken up time at yard $k$, in hour. |
| $t_k^{\text{classified}}$ | The parameter of the train classified time at yard $k$, in hour. |
| $t_k^{\text{delay}}$ | The parameter of the intermediate delay time at yard $k$, in hour. |
| $c_k^{\text{classified}}$ | The parameter of the train classified costs at yard $k$, in yuan. |
| $P_1^n$ | The base price 1 for different categories $n$ of delivering commodities (accumulation handling costs at the terminal without distance), in yuan/ton. |
| $R_2^n$ | The base price 2 for different categories $n$ of delivering commodities (shipping costs within distance), in yuan/ton-kilometer. |
| $L$ | The shipping distance, in kilometer. |
| $D$ | The freight demand during the contract's period, in ton. |
| $\Delta R$ | The difference between discounts provided the railroad and inventory costs of the customer, in yuan. |
| $C^{inventory}$ | Total inventory costs of the customer during the contract's period, in yuan. |
| $c^{inventory}$ | The unit inventory costs of the customer, in yuan/ton-day. |
| $T$ | The length of the contract period, in day. |

4.2 Formulation

$\Delta E$ can also be expressed by a series of differences in loading, unloading, reclassification and car-miles cost whether to operate the low-frequency entire train during the contract's period. We usually choose Eq.7 to calculate it.

$$\Delta E = \Delta E^{loading} + \Delta E^{unloading} + \Delta E^{reclassification} + \Delta E^{car-miles} \qquad (6)$$

When $\Delta E \geq 0$, it means that it is advantageous for the railroad to provide the low-frequency entire train for customers, since the total costs of operating the low-frequency entire train are less than the total consumption costs of operating the non-entire train. Otherwise, the railroad will not provide the low-frequency entire train. Therefore, $\Delta E$ is a prerequisite and the positive and negative value of $\Delta E$ represents whether railroads are willing to provide the low-frequency entire train.

Due to the car-miles cost is mostly equal comparing the low-frequency entire train service and the non-entire transportation, thus $\Delta E^{car-miles}$ can be ignored. The remaining parameters can be defined as follows.

$$\Delta E^{loading} = \gamma \Delta G^{loading} - C^{loading} T N_{ij} \qquad (7)$$

Comparing with the non-entire train service, operating the low-frequency entire train will increase extra loading costs. In addition, since most empty cars are supplied in groups at the loading station, the difference in loading time is not taken into account in the railroad practices, thus $\Delta G^{loading}$ can be ignored. Therefore, $\Delta E^{loading} < 0$.

$$\Delta E^{unloading} = \gamma \Delta G^{unloading} - C^{unloading} T N_{ij} \qquad (8)$$

Being similar to the situation at the loading station, when the number of unloading railcars is compatible with the unloading capacity, we can consider that the difference in unloading time is not taken into account, thus $\Delta G_{unloading}$ can also be ignored. Similarly, operating the low-frequency entire train will increase extra unloading costs, and $\Delta E^{unloading} < 0$.

Let $\Delta G^{reclassification}$ be the saving reclassification time.

$$t_k^{delay} = t_k^{broken\ up} + t_k^{classified} \qquad (9)$$

$$\Delta G^{reclassification} = TN_{ij} \sum_{k=1}^{q} t_k^{delay} \qquad i,j \in V \qquad (10)$$

Let $\Delta C^{reclassification}$ be the reducing intermediate handling costs.

$$\Delta C^{reclassification} = TN_{ij} \sum_{k=1}^{q} c_k^{classified} \qquad i,j \in V \qquad (11)$$

The low-frequency entire trains have been provided between the loading and unloading station, without reclassification at the classification yards. It is a key part for the railroad to reduce reclassification costs and save intermediate delay time, then we have:

$$\Delta E^{reclassification} = \gamma \Delta G^{reclassification} + \Delta C^{reclassification} \tag{12}$$

Customers calculate rail charges according to the following formula in China. The shipping distance $L$ refers to the shortest path from its origin to its destination on the line of China Railway network. The value of $P_1^n$ and $R_2^n$ for the different categories $n$ are derived according to the relative rules.

$$P = (P_1^n + R_2^n * L) * D \tag{13}$$

$C^{inventory}$ is calculated based on the following formula, and $\Delta S$ is calculated according to Section 3.2.

$$C^{inventory} = c^{inventory} T \Delta S \tag{14}$$

$\Delta H$ can be defined by:

$$\Delta H = \Delta E - \beta P \tag{15}$$

Where $\beta F$ denotes the amount of discount to encourage customers to choose the low-frequency entire train. When the value is positive, it means that the residual operating revenue except given discount is sufficient to cover the railroad's operating costs. Otherwise, the low-frequency entire train service with pricing incentive is not provided. The positive and negative value of $\Delta H$ represents whether railroads find a suitable discount to entice customers choosing the low-frequency entire train.

$\Delta R$ can be stated as:

$$\Delta R = \beta P - C^{inventory} \tag{16}$$

When $\Delta R \geq 0$, it means that amount of discount could compensate inventory costs of customers and customers will choose the low-frequency entire train service, conversely, the customer will not choose the low-frequency entire train service. The positive and negative value of $R$ represents whether customers are likely to choose the low-frequency entire train.

4.3 The trade-off strategies

There are several cases following the above analysis.

(1) With $\Delta E < 0$, it is clear that $\Delta H < 0$, in which the railroad still chooses to provide non-entire trains during the contract's period;

(2) With $\Delta E \geq 0$, if $\beta$ takes a small value, the situation $\Delta H \geq 0$ and $\Delta R < 0$ arises, in which the railroad expects the customer to choose the entire train service during the contract's period, but it is difficult for the customer to choose the entire train service.

(3) With $\Delta E \geq 0$, if $\beta$ takes a large value, it results in a situation where $\Delta H < 0$ and $\Delta R \geq 0$. In this state, the customer wants to choose the entire train service during the contract's period, and it is obvious that the railroad is not willing to provide the entire train service.

(4) With $\Delta E \geq 0$, if $\beta$ takes an appropriate value, only then will the situation of $\Delta H \geq 0$ and $\Delta R \geq 0$ arisen. This state is an ideal win-win situation, in which railroads are willing to provide the entire train service and customers are willing to choose the entire train service. In this case, we aid railroads in determining a suitable discount pricing parameter $\beta$ considering interests with customers.

Overall, this paper will determine reasonable discounts for operating the low-frequency entire train services, and measure the proportion of customers who are willing to choose the entire train services based on different discounts. Eventually, we provide pricing policy recommendations for operating the low-frequency entire train services.

## 5. Conclusion

This paper demonstrates the feasibility of operating low-frequency entire trains by analyzing the operating costs saving of railroads and the inventory costs increasing of customers. The railroad implements the optimal pricing scheme with a rail charges discount to maximize the proportion of the entire train. At the same time, each customer chooses its available train services to minimize the overall costs. The trade-off strategy is applied to keep a balance between the railroads and customers. Finally, the effectiveness of the entire train services adopting the discount pricing is examined. The paper shows the rail charges discount pricing strategy could increase the proportion of the entire train services without making customers worse off. The results and insights obtained provide valuable guidelines for railroads.

This study innovatively proposes a pricing strategy for China Railway and provides new ideas on rail charges rules. Furthermore, the discount pricing scheme could be further extended to allow for transportation demand uncertainties. The competition and cooperation of other transportation modes for freight and the influence of the pricing incentive strategy proposed in the paper deserve further research.

**Acknowledgments**: This study is supported by the National Key R&D Program of China (2018YFB1201402).

**Conflicts of Interest**: The authors declare no conflict of interest.

transportation, Operations research 62(2),383-400.